\documentclass{article}
\usepackage{amsmath,amsthm,amsfonts}
\usepackage{amssymb}
\usepackage{epsfig}
\def\Imp{\mathop{Im}}

\def\e{\epsilon}

\def\C{{\mathbb C}}

\def\CH{{\cal H}}
\def\CN{{\cal N}}
\def\CU{{\cal U}}

\def\N{{\mathbb N}}
\def\R{{\mathbb R}}

\def\tr{\mathop{tr}}

\def\diag{\mathop{diag}}

\renewcommand{\N}{{\bf N}}
\newcommand{\ignor}[1]{}
\newtheorem{lemma}{Lemma}
\newtheorem{claim}{Claim}
\newtheorem{theorem}{Theorem}
\newtheorem{problem}{Problem}
\newtheorem{remark}{Remark}
\newtheorem{corollary}{Corollary}
\newtheorem{definition}{Definition}
\begin{document}
\title{Almost commuting matrices with respect
to normalized Hilbert-Schmidt norm.}
\author{Lev Glebsky}
\maketitle
\begin{abstract}
Almost-commuting matrices with respect to the normalized
Hilbert-Schmidt norm are considered. Normal almost commuting
matrices are proved to be near commuting.
\end{abstract}
\section{Introduction}

Let $\C_{n\times n}$ be the set of complex $n\times n$-matrices. Let
$\CH_n,\CU_n,\CN_n\subset \C_{n\times n}$ be the sets of self-adjoin (Hermitian),
unitary and normal matrices, correspondingly. For $X,Y\in\C_{n\times
n}$ let $[X,Y]=XY-YX$. The following problem is classic
\begin{problem}\label{main}\textbf{(Must almost commuting matrices be nearly commuting?)}\\
Let $S_n=\C_{n\times n}, \CH_n, \CU_n$ or $\CN_n$.\\
For each $\e>0$, is there a $\delta=\delta(\e)>0$ such that for each
positive integer $n$, if $A,B\in S_n$ with
$\|A\|,\|B\|\leq 1$ and $\|AB-BA\|\leq\delta$, then there exist
$\tilde A,\tilde B\in S_n$ with $[\tilde A,\tilde B]=0$ and
$\|A-\tilde A\|,\|B-\tilde B\|\leq\e$?
\end{problem}
Here $\delta=\delta(\e)$ is independent of $n$; the non-uniform
version of the problem ($\delta=\delta(n,\e)$) has affirmative
answers \cite{Bast,Lux,Pear}. In spite of the equivalence of norms for finite
dimensional spaces, the answer on the uniform problem depends on the
norms $\|\cdot\|_n:\C_{n\times n}\to \R$. Indeed, the equivalence
may be non-uniform with respect to $n$. In the series of works
\cite{Choi1, Hastings,Lin,Voicul,Szarek}
the complete answer on Problem~\ref{main} have been found for
$\|\cdot\|=\|\cdot\|_{op}$. Where
$$
\|A\|_{op}=\sup\{\|Ax\|\;:\;\|x\|=1\}.\footnote{We consider $\C^n$ as a  Hilbert space
with scalar product $(x,y)=\sum x^*_iy_i$. It defines the Hilbert norm on $\C^n$,
$\|x\|=\sqrt{(x,x)}$. }
$$
The answer on Problem~\ref{main} is affirmative for $S_n=\CH_n$ and
negative for all other cases ($\|\cdot\|=\|\cdot\|_{op}$). We didn't
know any results for other norms. In the present paper we consider
Problem~\ref{main} for $\|\cdot\|=\|\cdot\|_{tr}$, the normalized
Hilbert-Schmidt norm:
$$
\|A\|_{\tr}=\sqrt{\frac{1}{n}\sum\limits_{j,k=1}^{n} |A_{j,k}|^2}.
$$
Our interest in normalized Hilbert-Schmidt norm arises from its use in
factor II
von Neumann algebras and hyperlinear groups, \cite{Neumann1, Pestov}.
It turns out that the normalized Hilbert-Schmidt norm is more
friendly for Problem~\ref{main}. We manage to prove that the answer
is affirmative for $S_n=\CH_n,\CU_n, \CN_n$, even if we speak about
several almost-commuting matrices. The reason why it is true is the
following. Transform the matrix $A$ into its diagonal form.
In this basis one
can approximate $\|\cdot\|_{tr}$-almost commuting matrices $A,B$ by
block diagonal matrices $\tilde A,\tilde B$ such that all blocks of
$\tilde A$ are multiples of unit matrices.

The technique of the paper is elementary. We systematically use that
the squire of  the normalized Hilbert-Schmidt norm of a block
diagonal matrix is a convex combination of the squires of the norms
of its blocks and concavity of some estimates, see
Section~\ref{sec_concave} for details.

All estimates in the theorems are given in the form
``$\|\cdot\|<C\e^\alpha$'', where $C$ is an integer. We do not try
to optimize the values of $C$, we just have decided that using some
proper numbers is less awkward than the use of expression of the type
``there exists $C>0$ such that...''.

\section{Notations and inequalities}
We consider $\C^n$ as a  Hilbert space with the scalar product
$(x,y)=\sum x^*_iy_i$. It defines the Hilbert norm on $\C^n$,
$\|x\|=\sqrt{(x,x)}$. Set $\C_{n\times n}$ of complex
$n\times n$ matrices naturally acts on $\C^n$.
As usual, we include $\C\subset\C_{n\times n}$ by constant diagonal matrices. So,
$1\in C_{n\times n}$, some times $1_n\in C_{n\times n}$ denotes the unit matrix.
For
$A=\{A_{i,j}\}\in\C_{n\times n}$ we define the  normalized trace
$$
\tr(A)=1/n\sum\limits_{i=1}^n A_{i,i}.
$$
It defines a scalar product on $\C_{n\times n}$:
$$
\langle A,B \rangle =\tr(A^*B)=\frac{1}{n}\sum_{i,j}A_{ij}^*B_{ij}
$$
and the normalized trace norm (normalized Hilbert-Schmidt norm)
$$
\|A\|_{tr}=\sqrt{\langle A,A \rangle}=\sqrt{\sum\limits_{i,j}|A_{ij}|^2}
$$
We also need the uniform operator norm
$$
\|A\|_{op}=\sup\{\|Ax\|\;:\;\|x\|=1\}
$$
We list some useful well-known inequalities, see \cite{Neumann1}, in the following:
\begin{lemma} \label{lm_ineq}
\begin{enumerate}
\item $|\langle A,B\rangle|\leq \|A\|_{tr}\|B\|_{tr}$
(the Cauchy-Schwarz inequality); substituting $1\to B$ gives
$|\tr(A)|\leq\|A\|_{tr}$;
\item $\|A+B\|_{tr}\leq \|A\|_{tr}+\|B\|_{tr}$
\item $\|AB\|_{tr}\leq \|A\|_{op}\|B\|_{tr}$ and
$\|BA\|_{tr}\leq \|A\|_{op}\|B\|_{tr}$
\item $\|A\|_{tr}\leq\|A\|_{op}\leq\sqrt{n}\|A\|_{tr}$
\item If $P$ is an orthogonal projector on k-dimensional subspace, then
$\|P\|_{tr}=\frac{\sqrt{k}}{\sqrt{n}}$.
\item If a matrix $A$ is of rank $k$, then there exists an orthogonal projector
$P$ of rank $k$ such that $A=PA$ and, by items 3,5
$\|A\|_{tr}\leq\sqrt{\frac{k}{n}}\|A\|_{op}$.
\item $\|A\|_{op}^2=\|A^*A\|_{op}$; $\|A\|_{tr}^2=\langle A,A\rangle=
\langle 1,A^*A\rangle\leq \|A^*A\|_{tr}$
\end{enumerate}
\end{lemma}
\begin{remark}
$\|\cdot\|_{op}$ is an algebraic norm: $\|AB\|_{op}\leq
\|A\|_{op}\|B\|_{op}$, but normalized trace norm is not a good
algebraic norm. We have only $\|AB\|_{tr}\leq
\sqrt{n}\|A\|_{tr}\|B\|_{tr}$ (it follows from 3,4 of Lemma~\ref{lm_ineq}).
\end{remark}
\begin{remark}
It is well known and easy to check that for unitary matrices $U,V$
$\|UXV\|_i=\|X\|_i$, where $i=op$ or $tr$. So, the norms $\|\cdot\|_i$ define
 norm on the set of linear operators from a Hilbert space to another
Hilbert space.
\end{remark}

The following lemma says that if $\|A\|_{tr}$ is
small then there is an orthogonal  projector  $P$ of a large rank,
such that $\|PA\|_{op}$ is small. Precisely,
\begin{lemma}\label{lemma_rank+norm}
For any $A\in C_{n\times n}$ there exists an orthogonal projector
$P$ such that
\begin{itemize}
\item $\|E-P\|_{tr}<\sqrt{\|A\|_{tr}},\;\;\;  \|PA\|_{op}<\sqrt{\|A\|_{tr}}$
\item If $A$ is normal, then, in addition, $AP=PA$.
\end{itemize}
\end{lemma}
\begin{proof}
Observe that $\|A\|_{op}= \sqrt{\|A^*A\|_{op}}=\sqrt{\|AA^*\|_{op}}$
and $AA^*$ is positive. Let
$0\leq\lambda_1\leq\lambda_2\leq...\leq\lambda_n$ be the eigenvalues
of $AA^*$. So,
$$
\|A\|_{tr}^2=\tr(AA^*)=\frac{1}{n}\sum_i \lambda_i \geq
\frac{\delta^2}{n}|\{i\;|\;\lambda_i\geq\delta^2\}|.
$$
So, we have $|\{i\;|\;\lambda_i\geq\delta^2\}|\leq
n\|A\|_{tr}^2/\delta^2$. Let $P_\delta$ be the orthogonal projector
on the space spanned by all eigenvectors of $AA^*$ with
$\lambda_i<\delta^2$. Then
$\|PA\|_{op}=\sqrt{\|PAA^*P\|_{op}}<\delta$ and $\|(E-P)\|_{tr}\leq
\|A\|_{tr}/\delta$. Putting $\delta=\sqrt{\|A\|_{tr}}$ proves the
first part of the lemma. The second part easily follows by
construction of $P$.
\end{proof}

\section{The concave estimate principle}\label{sec_concave}
We will need the following
\begin{claim}\label{claim_concave}
Let $\phi:\R^+\to\R^+$ be a concave increasing function with $\phi(0)=0$.
Then $x\to \phi^2(\sqrt{x})$ is a concave increasing function.
\end{claim}
\begin{proof}
On the set of functions $\R^+\to R^+$ we define an operation
$T$: $(T\phi)(x)=\phi^2(\sqrt{x})$. It is clear that if $\phi_1\leq
\phi_2$ then $T\phi_1\leq T\phi_2$. Here $\phi_1\leq \phi_2$ if
$\phi_1(x)\leq\phi_2(x)$ for all $x\in\R^+$. For all $x_0\in\R^+$
there exists $\alpha,\beta\geq 0$ such that $\phi(x_0)=\alpha
x_0+\beta$ and $\phi(x)\leq \alpha x+\beta$. (We have used here that $\phi(0)=0$.)
Observe that
$T(\alpha
x+\beta)=\alpha^2 x+\beta^2+2\alpha\beta\sqrt{x}$ is concave. So, for
any $x\in\R^+$ there exists a concave $f_x$, such that
$T\phi(x)=f_x(x)$ and $T\phi\leq f_x$. We deduce that $T\phi$ is
concave.
\end{proof}
Let $P\in\C[x_1,x_2,\dots,x_k,x_1^*,x_2^*,\dots,x_k^*]$ (polynomial
with complex coefficients, where $x_i^*$ is interpreted as complex
conjugate to $x_i$).
\begin{definition}\label{def1}
\begin{enumerate}
\item We say that matrices $A_1,..,A_k$ are an $\e$-solution of $P$
($\e$-satisfy $P$) if
$$
\|P(A_1,...,A_k,A_1^*,\dots,A_k^*)\|_{tr}\leq\e
$$
\item Let $\delta:\R^+\to \R^+$ be such that $\lim\limits_{\e\to 0}\delta_\e=0$.
Polynomial $P$ is called $\delta_\e$-stable if for any $\e$-solution
$A_1,A_2,...,A_k\in \C_{n\times n}$ of $P$ there exists an exact
solution $\tilde A_1,\tilde A_2,...,\tilde A_k\in \C_{n\times n}$ of
$P$ $\delta_\e$-close to $A_1,\dots A_k$, that is $\|A_i-\tilde
A_i\|_{tr}<\delta_\e$ for $i=1,...,k$. Polynomial $P$ is called stable if
it is $\delta_\e$-stable for some $\delta_\e $, $\lim\limits_{\e\to
0}\delta_\e=0$. (Note that $\delta_\e$ is independent of  $n$.)
\end{enumerate}
\end{definition}
Let ${\cal I}=\langle I_1,I_2,...,I_r\rangle$ be an ordered partition of $\{1,...,n\}$.
A matrix $A$ is said to be $\cal I$-{\bf block diagonal} (or just
$\cal I$-{\bf diagonal}) if nonzero elements of $A$ appear on
$I_j\times I_j$ places only. (It is clear that this is a usual block-diagonal matrix,
after conjugation with a permutation.) Similarly, we
call matrix $A$ {\bf cyclically $\cal I$-three-diagonal} if nonzero elements of
$A$ appear on $I_{j\oplus 1}\times I_j$, $I_j\times I_j$ and
$I_j\times I_{j\oplus 1}$ places only. Here $\oplus$ is the sum
$\mod r$.

Solution ($\e$-solution) $A_1,...,A_k$ of a polynomial $P$ is called
$\cal I$-diagonal solution ($\cal I$-diagonal $\e$-solution), if all
matrices $A_1,...,A_k$ are $\cal I$-diagonal. Under our assumptions,
$P(A_1,\dots,A_k^*)$ is $\cal I$-diagonal if all $A_i$ are.
\begin{lemma} \label{lm_princ_conv}
Suppose, that a polynomial $P$ is $\delta(\e)$-stable with a concave
$\delta(\e)$. Then for any $\cal I$-diagonal $\e$-solution of $P$
there exist an $\cal I$-diagonal solution of $P$ that is
$\delta(\e)$-close to this $\e$-solution.
\end{lemma}
\begin{proof}
The proof uses the following  facts:
\begin{enumerate}
\item $\sum_j \alpha_j\delta(x_j)\leq
\delta (\sum_j\alpha_jx_j)$ for $x_j,\alpha_j\geq 0$, $\sum_j\alpha_j=1$ and
concave $\delta$.
\item Let $A$ be an $\cal I$-diagonal matrix and $A^1,A^2,\cdots,A^r$ its
diagonal components of dimensions $d_1,\dots,d_r$, correspondingly.
Then
$$
\|A\|_{tr}^2=\sum_j\frac{d_j}{n}\|A^j\|_{d_j}^2,
$$
where $n=d_1+d_2+\cdots+d_r$ and $\|\cdot\|_{d}$ is the normalized
Hilbert-Schmidt norm on $\C_{d\times d}$.
\end{enumerate}
Let $A_1,\dots,A_k$ be an $\cal I$-diagonal solution of $P$ with
diagonal components $A_j^l$. Let $\|P(A_1^l,...,A_k^l)\|=\e_l$. We
have $\epsilon^2\geq\sum\frac{d_l}{n}\e_l^2$. There exists a
solution $\tilde A_1^l,\dots,\tilde A_k^l$ of $P$ with
$$
\|\tilde A_j^l-A_j^l\|_{tr}^2\leq\delta^2(\sqrt{\e^2}).
$$
Solution $\tilde A_j$ is constructed by blocks $\tilde A_j^l$. Now
$$
\|\tilde A_j-A_j\|^2_{tr}=\sum_l\frac{d_l}{n}\|\tilde
A_j^l-A_j^l\|^2\leq \sum_l \frac{d_l}{n}\delta^2(\sqrt{\e_l^2})\leq
\delta^2\left(\sqrt{\sum_l\frac{d_l}{n}\e_l^2}\right)\leq \delta^2(\e).
$$
Here we use concavity of $\delta^2(\sqrt{x})$ by
Claim~\ref{claim_concave}.
\end{proof}

\section{Almost unitary matrices are near unitary}

\begin{lemma}\label{lemma_almost_unitary_op}
Let $B:L_1\to L_2$ be an unitary operator from a
Hilbert space $L_1$ to a Hilbert space $L_2$ such
that $\|B^*B-1_{L_1}\|_{op}\leq\epsilon\leq 1/3$. Then there exists
an unitary operator $V:L_1\to L_2$ such that $\|B-V\|_{op}< 2\epsilon$.
\end{lemma}

\begin{proof}
Just take $V=B(B^*B)^{-1/2}$, where
$(B^*B)^{-1/2}=(1_{L_1}-X)^{-1/2}=
\sum\limits_{j=0}^\infty \frac{(2j)!X^j}{2^{2j}(j!)^2}$.
(We have denoted $X=B^*B-1$.)
Make the following estimates
$$
\|(B^*B)^{-1/2}-1\|_{op}\leq
\sum\limits_{j=1}^\infty \frac{(2j)!\epsilon^j}{2^{2j}(j!)^2}<
\sum\limits_{j=1}^\infty \epsilon^j=\frac{\epsilon}{1-\epsilon}.
$$
Now,
$$
\|V-B \|_{op}\leq \|B\|_{op}\|(B^*B)^{-1/2}-1\|_{op}<
\frac{\epsilon(\epsilon+1)}{1-\epsilon}\leq 2\epsilon,
$$
for $\epsilon\leq 1/3$. (We were using the fact that
$\|B\|_{op}=\sqrt{\|B^*B\|_{op}}\leq\sqrt{\|B^*B-1\|_{op}+1}\leq\|B^*B-1\|_{op}+1$.)
\end{proof}

\begin{lemma}\label{lemma_almost_unitary_tr}
Let $\|A^*A-1\|_{tr}\leq\epsilon\leq 1/3$ for a matrix $A$.
Then there exists a unitary $U$, such that
$\|A-U\|_{tr}\leq 5\epsilon^\frac{1}{4}$ and $\|A-U\|_{tr}\leq
(3+\|A\|_{op})\sqrt{\epsilon}$.
\end{lemma}

\begin{proof}
By Lemma~\ref{lemma_rank+norm} there exists orthogonal projector
$P$, $\|1-P\|_{tr}\leq\sqrt{\epsilon}$ such that
$\|PA^*AP-P\|_{op}\leq\sqrt{\epsilon}$. Let $X=\Imp P$. Consider the
restriction $B=A|_X :X\to Y=A(X)$, then $B^*:Y\to X$. Observe
that $B^*=PA^*|_Y$. So, $\|B^*B-1_X\|_{op}\leq\sqrt{\epsilon}$
and, by Lemma~\ref{lemma_almost_unitary_op}, there exists a unitary
$V:X\to Y$ with $\|V-B\|_{op}<2\sqrt{\epsilon}$. Let $\tilde V$
be any unitary operator from $X^\perp$ to $Y^\perp$. Take
$U=V\oplus\tilde V$. We estimate:
$$
\|A-V\oplus\tilde V\|_{tr}\leq \|P_YAP_X-P_YVP_X\|_{tr}+
\|AP_X^\perp\|_{tr}+\|P_Y^\perp\tilde VP_X^\perp\|_{tr}\leq
$$
$$
\|P_YAP_X-P_YVP_X\|_{op} +\|\tilde V\|_{op}\|P_X^\perp\|_{tr}
+\|AP_X^\perp\|_{tr}\leq 3\sqrt{\epsilon}+\|AP_X^\perp\|_{tr}.
$$
Now, the first inequality of the lemma follows from
$$
\|AP_X^\perp\|_{tr}^2\leq\|P_X^\perp A^*AP_X^\perp\|_{tr}\leq
\|P_X^\perp\|_{tr}+ \|P_X^\perp A^*AP_X^\perp-P_X^\perp\|_{tr}\leq
$$
$$
\sqrt{\epsilon}+
\|P_X\|_{op}^2\|A^*A-1\|_{tr}\leq \sqrt{\epsilon}+\epsilon.
$$
The second inequality of the lemma follows from
$$
\|AP_X^\perp\|_{tr}\leq \|A\|_{op}\|P_X^\perp\|_{tr}\leq \|A\|_{op}\sqrt{\epsilon}.
$$
\end{proof}
We will work with $\cal I$-diagonal matrices, so we
need a global concave estimate.
\begin{corollary}\label{cor_unit}
Let $\|A\|_{op}\leq 3$. Then there exists a unitary matrix $V$ such
that $\|A-V\|_{tr}\leq 6\sqrt{\|A^*A-1\|_{tr}}$.
\end{corollary}
\begin{proof}
For $\|A^*A-1\|_{tr}\leq 1/3$ it is
Lemma~\ref{lemma_almost_unitary_tr}. Further, $\|A-V\|_{tr}\leq
\|A\|_{tr}+1$ and $\|A\|_{tr}^2\leq \|A^*A\|_{tr}\leq
\|A^*A-1\|_{tr}+1$. So, $\|A-V\|_{tr}\leq\sqrt{\|A^*A-1\|_{tr}+1}+1$.
It remains to check
that $6\sqrt{x}>\sqrt{x+1}+1$ for $x>1/3$.
\end{proof}
\section{Almost commuting unitary matrices are near commuting}
\begin{theorem}\label{th_2_unit}
Let $U_1$ and $U_2$ be unitary matrices. Then there exists unitary matrices
$A_1,A_2$, $[A_1,A_2]=0$ such that $\|U_1-A_1\|_{tr}\leq 30\left(\|[U_1,U_2]\|_{tr}\right)^{1/9}$
and $\|U_2-A_2\|_{tr}\leq 30\left(\|[U_1,U_2]\|_{tr}\right)^{1/9}$. In addition, $[A_1,U_1]=0$.
\end{theorem}
Before the proof of the theorem we consider an example where $U_2$
is a cyclic permutation and $U_1$ its diagonal form:
$U_1=\diag(w,w^2,\dots,w^n=1)$ with $w=\exp(\frac{2\pi i}{n})$ and
$U_2=P_n$ with
$$
P_{j,k}=\left\{\begin{array}{lll} 1 &\mbox{if} & j=k+1\mod n \\
                                  0 &\mbox{overwise}
               \end{array}\right..
$$
This is a counterexample  to Problem~\ref{main}  for
$\|\cdot\|=\|\cdot\|_{op}$ and $S_n=\CU_n$ found by Voiculescu,
\cite{Voicul}
(he proves that it is indeed a counterexample). One has
$\|[U_1,U_2]\|_{op}=\|[U_1,U_2]\|_{tr}=|1-w|\to 0$ when
$n\to\infty$. Suppose, for simplicity, that $n=md$ for large $m$ and
$d$. Then we can take as $A_1$ and $A_2$ the following
block-diagonal matrices: $A_1=\diag(\tilde w1_d,\tilde w^21_d,\dots,
\tilde w^m1_d)$, where $\tilde w=\exp(\frac{2\pi i}{m})$ and
$A_2=\diag(P_d,P_d,\dots,P_d)$.
\begin{proof}
After transformation $U_1\to V^{-1}U_1V$, $U_2\to V^{-1}U_2V$, assume
that $U_1=\diag(\alpha_1,\alpha_2,....,\alpha_n)$. The main idea of
the proof is the following. Change some elements of $U_2$ by $0$ and
approximate $U_1$ by a diagonal matrix with spectrum $\exp(2\pi
i\frac{j}{m})$ for proper $m$ in such a way that $U_1$ and $U_2$
become block diagonal matrices with all blocks of $U_1$ being
multiples of unit matrices. New $U_1$ and $U_2$ are commuting, but now
$U_2$ is not unitary. Approximate $U_2$ by an unitary matrix,
conserving its block structure. It can be done using
Corollary~\ref{cor_unit} and Lemma~\ref{lm_princ_conv}. Let us
describe this procedure in details.

Let $\|[U_1,U_2]\|_{tr}=\e$.
 Take a positive
integer $t\geq 6$ that will be optimized latter. Let
$w=\exp(\frac{2\pi i}{t})$. Let $|1-w|=\Delta$. One has
$$
6/t\leq\Delta\leq 2\pi/t<7/t.
$$
\begin{enumerate}
\item Let $U_2=\{u_{jk}\}$. Define $\tilde U_2=\{\tilde u_{jk}\}$ by
the following rule:
$$
\tilde u_{jk}=\left\{\begin{array}{lll}
             u_{jk} &\mbox{if} & |\alpha_j-\alpha_k|<\Delta\\
             0     &\mbox{if} & |\alpha_j-\alpha_k|\geq\Delta
              \end{array}\right.
$$
One has
$$
\|\tilde U_2-U_2\|_{tr}\leq\e/\Delta
$$
Indeed,
$$
n\e^2\geq n\|[U_1,U_2]\|^2=\sum\limits_{j=1,k=1}^{n,n}
|\alpha_k-\alpha_j|^2|u_{jk}|^2\geq
\Delta^2\sum\limits_{|\alpha_k-\alpha_j|\geq\Delta} |u_{jk}|^2=\Delta^2 n\|\tilde
U_2-U_2\|_{tr}^2
$$
One can check that $\|[U_1,\tilde U_2]\|\leq \|[U_1,U_2]\|\leq\e$.
\item  Approximate $U_1$ by a diagonal
matrix $\tilde U_1$ with spectrum $\{w^j\;:\;j=0,1,\dots,t-1\}$.
Precisely, let ${\cal I}=\{I_0,I_1,\dots,I_{t-1}\}$ with
$I_j=\{l\;:\;\alpha_l\in (w^{j-1/2},w^{j+1/2}]\}$, where $(x,y]$ is a
semiopen arc of the unit circle in $\C$. Then $\tilde U_1$ is an $\cal
I$-diagonal matrix with $j$ block $U^j_1=w^j$. One has that
$$
\|\tilde U_1-U_1\|_{tr}\leq\Delta
$$
and $\|[\tilde U_1,\tilde U_2]\|_{tr}\leq \e+2\Delta$. Observe that
$\tilde U_2$ is a cyclically $\cal I$-three-diagonal matrix.
\item Fix another parameter $a\in\N$ that will be optimized latter.
One can find $S=\{s_0,s_2,\dots,s_{c-1}\}\subset\{0,1,2,\dots,t-1\}$ such that
\begin{itemize}
\item $a\leq |s_{r\oplus1}-s_r|\leq 3a$, for any $r=0,1,...,c-1$.
Here $\oplus$ is the sum $\mod c$.
\item $|I_j|\leq \frac{n}{a}$ for any $j\in S$.
\item $|S|=c\leq t/a$
\end{itemize}
\item In order to construct  $A_1$ we make a more rough partition
${\cal\tilde I}=\{\tilde I_0,\tilde I_2,\dots,\tilde I_{c-1}\}$.
Where $\tilde I_j=I_{s_j}\cup I_{s_j+1}\cup\cdots\cup I_{s_{j\oplus
1}-1}$, where $\pm$ is $\mod t$ and $\oplus$ is $\mod c$. Now, $A_1$ is
a cyclically $\cal\tilde I$-diagonal matrix with the blocks
$A_1^j=w^{\frac{1}{2}(s_j+s_{j\oplus}-1)}$.
By the first item of 3) we have
$$
\|\tilde U_1-A_1\|_{tr}\leq |1-w^{\frac{3}{2}a}|\leq \frac{3}{2}a\Delta
$$
and $\|[A_1,\tilde U_2]\|_{tr}\leq \e+2\Delta+3a\Delta$.
\item We give our construction of $A_2$ in two steps.
Recall that $\tilde U_2$ is $\cal I$-three-diagonal. We construct $B$ by
removing  from $\tilde U_2$ the blocks $I_{j-1}\times I_j$ and
$I_j\times I_{j-1}$ for each $j\in S$. The resulting matrix $B$ is
$\cal\tilde I$-diagonal and, consequently, $[A_1,B]=0$. We estimate:
$$
\|\tilde U_2-B\|_{tr}=\sum_{j\in S}
(\|U_2^{j-1,j}\|_{tr}+\|U_2^{j,j-1}\|_{tr}) \leq
2\frac{|S|}{\sqrt{a}}\leq 2\frac{t}{a^{3/2}}.
$$
For the first inequality we use
$\|U_2^{j-1,j}\|_{tr}\leq\sqrt{\frac{|I_j|}{n}}\|U_2^{j-1,j}\|_{op}$
(the item 6. of Lemma~\ref{lm_ineq}) and $\|U_2^{j-1,j}\|_{op}\leq
\|U_2\|_{op}=1$ (the operator norm of a submatrix is less than the
operator norm of the matrix). The same inequalities are valid for
$\|U_2^{j,j-1}\|$. The second inequality is a property of $S$.
\item The matrix $B, A_1$ are $\cal\tilde I$-diagonal and each block of
$A_1$ is
a multiple of the unit matrix, so $[A_1,B]=0$. The problem is that
$B$ is not unitary.
$\|U_2-B\|_{tr}\leq\frac{\e}{\Delta}+2\frac{t}{a^{3/2}}\leq
\e t/6+2ta^{-3/2}=\gamma$.
and $\|B\|_{op}\leq 3$ ($B$ is $\cal I$-three-diagonal with the
operator norm of each block less than $1$ as submatrices of a
unitary matrix.) It follows that
$$
\|B^*B-1\|_{tr}=\|B^*B-U_2^*U_2\|_{tr}\leq \|B^*B-B^*U_2\|_{tr}
+\|B^*U_2-U_2^*U_2\|_{tr}
$$
$$
\leq \|B^*\|_{op}\|B-U_2\|_{tr}+\|U_2\|_{op}\|B^*-U_2^*\|_{tr}\leq
4\gamma.
$$
The matrix $B$ is an $\cal\tilde I$-diagonal matrix. By
Lemma~\ref{lm_princ_conv} and Corollary~\ref{cor_unit} there exists
a unitary $\tilde I$-diagonal matrix $A_2$ with
$$
\|B-A_2\|_{tr}\leq 12\sqrt{\gamma}=12\sqrt{\frac{\e
t}{6}+2ta^{-3/2}}
$$
It is clear that $[A_1,A_2]=0$.
\item We only need to choose $a$, $t$ and estimate $\|U_i-A_i\|_{tr}$.
Suppose, for a moment\footnote{The condition on $\e$
is to guarantee $t\geq 6$.}, that $\e\leq 6^{-9/7}$, choose $a,t\in\N$ such that
$\e^{-7/9}\leq t\leq 2\e^{-7/9}$ and $\e^{-2/3}\leq a\leq
2\e^{-2/3}$. We have:
$$
\|U_1-A_1\|_{tr}\leq \Delta+\frac{3}{2}a\Delta\leq
\frac{7}{t}+\frac{21a}{2t}\leq 7\e^{7/9}+21\e^{1/9}<30\e^{1/9}.
$$
Further,
$$
\|U_2-A_2\|_{tr}\leq \frac{1}{6}\e t+2ta^{-3/2}+12
\sqrt{\frac{1}{6}\e t+2ta^{-3/2}}\leq \frac{13}{3}\e^{2/9}+
12\sqrt{\frac{13}{3}}\e^{1/9}< 30\e^{1/9} .
$$
For $\e\geq 6^{-9/7}$ we have
$$
\|A_i-U_i\|_{tr}\leq 2\leq 30\cdot6^{-1/7}\leq 30\e^{1/9}
$$
\end{enumerate}
The pair $A_1,A_2$ satisfies the statement of the theorem.
\end{proof}
\ignor{In fact we prove
\begin{corollary}
Let $U_1$ and $U_2$ be unitary matrices, such that
$\|[U_1,U_2]\|_{tr}\leq\e$. Then there exist a partition $\cal I$ and unitary matrix $V$
such that $\|V^{-1}U_1V-A_1\|_{tr}\leq 26\e^{1/9}$
and $\|V^{-1}U_2V-A_2\|_{tr}\leq 26\e^{1/9}$ for $\cal I$-block diagonal
unitary $A_1, A_2$ with
all blocks of $A_1$ being multiples of  unit matrices.
\end{corollary}
Concavity of $x^{1/9}$ with Lemma~\ref{lm_princ_conv} gives
\begin{corollary}
Let $U_1$ and $U_2$ be $\cal I$-block unitary matrices, such that
$\|[U_1,U_2]\|_{tr}\leq\e$. Then there exists $\cal I$-block unitary
matrix $V$ and a partition $\cal I'$, a refinement of $\cal I$, such
that $\|V^{-1}U_1V-A_1\|_{tr}\leq 26\e^{1/9}$ and
$\|V^{-1}U_2V-A_2\|_{tr}\leq 26\e^{1/9}$ for $\cal I'$ diagonal
unitary $A_1$ and $A_2$ with all blocks of $A_1$ being multiple of
unit matrices.
\end{corollary}
}
We need the following
\begin{claim}\label{Cl_commut_ineq}
Let $\|A\|_{op},\|B\|_{op},\|\tilde A\|_{op},\|\tilde B\|_{op}\leq 1$.
Then
$$
\|[\tilde A,\tilde B]\|_{tr}\leq \|[A,B]\|_{tr}+
2(\|A-\tilde A\|_{tr}+\|B-\tilde B\|_{tr}).
$$
\end{claim}
\begin{proof}
$\|AB-\tilde A\tilde B\|_{tr}=\|AB-A\tilde B +A\tilde B-\tilde A\tilde B\|_{tr}\leq
\|A\|_{op}\|B-\tilde B\|_{tr}+\|A-\tilde A\|_{tr}\|\tilde B\|_{op}\leq
\|B-\tilde B\|_{tr}+\|A-\tilde A\|_{tr}.$ Combining it with the same estimate
for $\|BA-\tilde B\tilde A\|_{tr}$ we get the claim.
\end{proof}
\begin{theorem}\label{th_unit}
There exists $\delta(\e,k)$, $\delta(\e,k)\to 0$ when $\e\to 0$ for any
$k\in\N$,
such that
if $\|[U_i,U_j]\|_{tr}\leq\e$ for unitary $U_1,U_2,\dots,U_k$,
then there exist pairwise commuting unitary matrices $A_1,\dots,A_k$
such that $\|U_j-A_j\|_{tr}\leq \delta(\e,k)$.
\end{theorem}
\begin{proof}

Let $\psi(x)=30x^{1/9}$ and $\phi_j(\cdot)$ be defined by the relation:
$$
\phi_0(x)=x,\;\;\;\;\phi_{j+1}(x)=4\psi(\phi_j(x))+x.
$$
For $r=1,\dots k-1$ we prove by induction the following statement
\\~\\~
{\it  There exist a unitary matrix $V$, a partition ${\cal I}_r$ of $\{1,\dots,n\}$, and
${\cal I}_r$-diagonal  matrices $\tilde
U_1,\tilde U_2,\dots,\tilde U_k$, such that
\begin{itemize}
\item All blocks of $\tilde U_1,\dots \tilde U_r$ are multiples of the unit
matrix.
\item $\|\tilde U_j-V^{-1}U_jV\|_{tr}\leq \psi\left(\phi_{j-1}(\e)\right)$,
for $j\leq r$ and $\|\tilde U_j-V^{-1}U_jV\|_{tr}\leq \psi\left(\phi_{r-1}(\e)\right)$,
for $j>r$.
\end{itemize}
}~\\
The theorem follows from the Statement for $r=k-1$ and
$\delta(\e,k)=\psi\left(\phi_{k-1}(\e)\right)$. Let us proof the
Statement.

{\bf $r=1$.} In the proof of Theorem~\ref{th_2_unit} matrix $A_1$
and partitions $\cal I$ and $\cal\tilde I$ is independent of $U_2$.
The construction of $A_2$ depends on partitions $\cal I$ and
$\cal\tilde I$ only. So, we may construct $\cal\tilde I$-diagonal
$\tilde U_1,\tilde U_2,\dots,\tilde U_k$ satisfying the Statement
for $r=1$.

{\bf $r\to r+1$.} Let $\tilde U_1,\dots, \tilde U_k$ be as in the
Statement. Then $\|[\tilde U_i,\tilde U_j]\|=0$ for $i< r$ and,
by Claim~\ref{Cl_commut_ineq},
$\|[\tilde U_i,\tilde U_j]\|\leq \phi_{r}(\e)$ for $i,j\geq
r$. Let $I_l\in {\cal I}_r$. Work with $\tilde U_r^l,\dots, \tilde
U_k^l$ as in the proof for $r=1$. Then apply
Lemma~\ref{lm_princ_conv}.
\end{proof}

\section{Self-adjoint matrices.}

For every almost-commuting self-adjoint matrices $A,B$ we construct
commuting self-adjoint matrices with the same operator norm and
close to $A,B$ by the normalized Hilbert-Schmidt norm. In order to
preserve the operator norm we need
\begin{lemma}\label{lm_self-adjoin_norm}
Let $A$, $B$, $C=A+B$ be self adjoint matrices. Let $D(B)$ and $D(C)$ be the
decreasing diagonal form of $B$ and $C$, correspondingly. Then
$\|D(C)-D(B)\|_{op}\leq \|A\|_{op}$
\end{lemma}
\begin{proof}
Let $\alpha_1\geq\alpha_2\geq\dots\geq\alpha_n$,
$\beta_1\geq\dots\geq\beta_n$ and $\gamma_1\geq\dots\geq\gamma_n$ be
(ordered) eigenvalues of A,B and C, correspondingly.  The H.Weyl
inequality \cite{Fulton,Weyl} states:
$$
\gamma_{j+k-1}\leq \alpha_{j}+\beta_{k}.
$$
Writing $-C=-A-B$ and reordering the eigenvalues we get:
$$
\gamma_{k-j+1}\geq \alpha_{n-j+1}+\beta_{k}.
$$
Putting $j=1$ in the both inequalities and using the fact that
$\alpha_1, -\alpha_n\leq \|A\|_{op}$ we get
$$
-\|A\|_{op}+\beta_k\leq\gamma_k\leq\|A\|_{op}+\beta_k.
$$
\end{proof}

\begin{corollary} \label{cor_norm}
Let $A,C$ be self-adjoint and $C$ be $\cal I$-diagonal. Then there
exists $\cal I$-diagonal self-adjoint matrix $\tilde C$ such that
$\|\tilde C\|_{op}\leq\|A\|_{op}$ and $\|\tilde C-C\|_{tr}\leq
\|C-A\|_{tr}$
\end{corollary}
\begin{proof}
We can choose $\tilde C$ such that $D(\tilde C)=D(C)-D(C-A)$.
\end{proof}
\begin{theorem}\label{th_2_self-adjoint}
Let $H_1$ and $H_2$ be self-adjoint matrices, such that
$\|H_i\|_{op}\leq 1$, $i=1,2$.
 Then there exists  self-adjoint matrices
$A_1,A_2$, $[A_1,A_2]=0$ such that
$\|H_1-A_1\|_{tr}\leq 12\left(\|[H_1,H_2]\|_{tr}\right)^{1/6}$
and $\|H_2-A_2\|_{tr}\leq 12\left(\|[H_1,H_2]\|_{tr}\right)^{1/6}$,
$\|A_i\|_{op}\leq 1$. In addition, $[A_1,H_1]=0$.
\end{theorem}
\begin{proof}
We follow the same routine as in the proof of
Theorem~\ref{th_2_unit}. Instead of
Lemma~\ref{lemma_almost_unitary_tr} we use Corollary~\ref{cor_norm} to
keep the operator norm.

Let $\|[H_1,H_2]\|_{tr}=\e$
We suppose that $H_1=\diag(\alpha_1,\alpha_2,....,\alpha_n)$,
$-1\leq\alpha_1\leq\alpha_2\leq\dots\leq\alpha_n\leq 1$.
Take a positive integer $t$ that will be optimized latter.
\begin{enumerate}
\item Let $H_2=\{h_{jk}\}$. Define $\tilde H_2=\{\tilde h_{jk}\}$ by
the following rule:
$$
\tilde h_{jk}=\left\{\begin{array}{lll}
             h_{jk} &\mbox{if} & |\alpha_j-\alpha_k|<1/t\\
             0     &\mbox{if} & |\alpha_j-\alpha_k|\geq 1/t
              \end{array}\right.
$$
As in the proof of Theorem~\ref{th_2_unit} one has
$$
\|\tilde H_2-H_2\|_{tr}\leq \e t.
$$
Clearly, $\tilde H_2$ is self-adjoin.
\item  Approximate $H_1$ by a diagonal
matrix $\tilde H_1$ with spectrum $\{\frac{j}{t}\;:\;j=-t,\dots,t\}$.
Precisely, let ${\cal I}=\{I_{-t},I_{-t+1},\dots,I_{t}\}$ with
$I_j=\{l\;:\;\alpha_l\in (\frac{2j-1}{2t},\frac{2j+1}{2t}]\}$, where $(x,y]$ is a
semiopen interval. Then $\tilde H_1$ is an $\cal
I$-diagonal matrix with the $j$-th block $H^j_1=\frac{j}{t}$. One has that
$$
\|\tilde H_1-H_1\|_{tr}\leq\frac{1}{t}
$$
and that $\tilde H_2$ is an $\cal I$-three-diagonal matrix (not cyclically
$\cal I$-three-diagonal).
\item Fix another parameter $a\in\N$ that will be optimized latter.
One can find $S=\{s_0,s_2,\dots,s_{c-1}\}$, $-t\leq s_1<s_1<\dots s_{c-1}\leq t$
such that
\begin{itemize}
\item $a\leq |s_{r+1}-s_r|\leq 2a$, for any $r=0,1,...,c-1$.
\item $|I_j|\leq \frac{n}{a}$ for any $j\in S$.
\item $|S|=c\leq (2t+1)/a$
\end{itemize}
\item In order to construct  $A_1$ we make more rough partition
${\cal\tilde I}=\{\tilde I_0,\tilde I_2,\dots,\tilde I_{c-1}\}$.
Where $\tilde I_j=I_{s_j}\cup I_{s_j+1}\cup\cdots\cup I_{s_{j+1}-1}$.
Now, $A_1$ is
cyclically $\cal\tilde I$-diagonal matrix with block
$A_1^j=\frac{s_j+s_{j+1}-1}{2t}$.
By the first item of 3) we have
$$
\|\tilde H_1-A_1\|_{tr}\leq \frac{a}{t}
$$
\item We give construction of $A_2$ in two steps.
Recall that $\tilde H_2$ $\cal I$-three-diagonal. We construct $B$ by
removing  from $\tilde H_2$ blocks $I_{j-1}\times I_j$ and
$I_j\times I_{j-1}$ for each $j\in S$. The resulting matrix $B$ is
$\cal\tilde I$-diagonal and, consequently, $[A_1,B]=0$. We estimate:
$$
\|\tilde H_2-B\|_{tr}=\sum_{j\in S}
(\|H_2^{j-1,j}\|_{tr}+\|H_2^{j,j-1}\|_{tr}) \leq
2\frac{|S|}{\sqrt{a}}\leq 2\frac{2t+1}{a^{3/2}}.
$$
For the first inequality we use
$\|H_2^{j-1,j}\|_{tr}\leq\sqrt{\frac{|I_j|}{n}}\|H_2^{j-1,j}\|_{op}$
(the item 6. of Lemma~\ref{lm_ineq})
and $\|H_2^{j-1,j}\|_{op}\leq \|H_2\|_{op}\leq 1$ (the operator norm of
a submatrix is less than the operator norm of the matrix). The same
inequalities are valid for $\|H_2^{j,j-1}\|$. The second inequality is
a property of $S$.
\item The matrix $B$ is $\cal\tilde I$-diagonal, self-adjoint, and
$$
\|H_2-B\|_{tr}\leq\e t + 2\frac{2t+1}{a^{3/2}}.
$$
So by Corollary~\ref{cor_norm} there exists $\cal\tilde I$-diagonal
self-adjoint $A_2$ with
$$
\|H_2-A_2\|_{tr}\leq 2\e t + 4\frac{2t+1}{a^{3/2}}.
$$
It is clear that $[A_1,A_2]=0$.
\item We only need to choose $a$, $t$ and estimate $\|H_i-A_i\|_{tr}$.
Suppose, for a moment, that $\e\leq 4^{-1}$, choose $a,t\in\N$ such that
$\frac{1}{2}\e^{-5/6}\leq t\leq \e^{-5/6}$ and $\e^{-2/3}\leq a\leq
2\e^{-2/3}$.
We have:
$$
\|H_1-A_1\|_{tr}\leq \frac{1}{t}+\frac{a}{t}\leq 2\e^{5/6}+4\e^{1/6}<
12\e^{1/6}
$$
Further,
$$
\|H_2-A_2\|_{tr}\leq 2\e^{1/6}+8\e^{1/6}+4\e<12\e^{1/6}
$$
For $\e\geq 4^{-1}$ we have
$$
\|A_i-H_i\|_{tr}\leq \|A_i-H_i\|_{op}\leq 2< 12(4^{-1/6})<12\e^{1/6}
$$
\end{enumerate}
The pair $A_1,A_2$ satisfies the statement of the theorem.
\end{proof}

\begin{theorem}\label{th_self-adjoint}
There exists $\delta(\e,k)$, $\delta(\e,k)\to 0$ when $\e\to 0$ for any $k\in\N$,
such that
if $\|[H_i,H_j]\|_{tr}\leq\e$ for self-adjoint matrices $H_1,H_2,\dots,H_k$ with
$\|H_i\|_{op}\leq 1$,
then there exist pairwise commuting self-adjoint matrices $A_1,\dots,A_k$
such that $\|U_j-A_j\|_{tr}\leq \delta(\e,k)$ and $\|A_i\|\leq 1$.
\end{theorem}
\begin{proof}
The same as for Theorem~\ref{th_unit}.
\end{proof}
\section{Normal matrices}

Observe that Theorem~\ref{th_self-adjoint} implies the existence of commuting
normal matrices
close to almost commuting ones.
Observe also, that Theorem~\ref{th_2_self-adjoint} implies the existence of a normal matrix
$N$ close to an $\|\cdot\|_{tr}$-almost normal matrix $M$.
Could it be done in a way that $\|N\|_{op}\leq\|M\|_{op}$?
In the section we give the affirmative answer to this question
(Corollary~\ref{cor_norm_cons})

\begin{theorem}\label{th_2_self-adjoint_and_unitary}
Let $U$ and $H$ be unitary and positive matrices, correspondingly.
Let $\|H\|_{op}\leq 1$.
 Then there exists  unitary and positive matrices
$V,A$ such that $[V,A]=[H,A]=0$, $\|V-U\|_{tr}\leq 30\|[U,H]\|_{tr}^{1/9}$
and $\|H-A\|_{tr}\leq 30\|[U,H]\|_{tr}^{1/9}$, $\|A\|_{op}\leq 1$.
\end{theorem}
\begin{proof}
Let $H=\diag(h_1,\dots,h_n)$. Make partition $\cal I$ and $\cal\tilde I$ as in the proof of
Theorem~\ref{th_2_self-adjoint}. Construct $A$ as $A_1$ in Theorem~\ref{th_2_self-adjoint}
and $V$ as $U_2$ in Theorem\ref{th_2_unit}.
\end{proof}
\begin{lemma}\label{lemma_root_estimate}
Let $A,B$ be positive commuting matrices. Then
$\|A-B\|_{tr}\leq\sqrt{\|A^2-B^2\|_{tr}}$.
\end{lemma}
\begin{proof}
Without loss of generality we may assume that $A=\diag(a_1,a_2,\dots,a_n)$
and $B=\diag(b_1,b_2,\dots,b_n)$. Now,
$$
\|A-B\|_{tr}^2=\frac{1}{n}\sum\limits_{j=1}^n(a_j-b_j)^2\mathop\leq\limits_{\bf (a)}
\frac{1}{n}\sum\limits_{j=1}^n|a_j^2-b_j^2|=
$$
$$
\frac{1}{n}\sum\limits_{j=1}^n\sqrt{(a_j^2-b_j^2)^2}\mathop\leq\limits_{\bf (b)}
\sqrt{\frac{1}{n}\sum\limits_{j=1}^n(a_j^2-b_j^2)^2}=\|A^2-B^2\|_{tr}
$$
The inequality {\bf (a)} is due to $(a-b)^2\leq |a^2-b^2|$ for $a,b\geq 0$; the inequality
{\bf (b)} is
due to concavity of $\sqrt{\cdot}$.
\end{proof}
Theorem~\ref{th_2_self-adjoint_and_unitary} with
Lemma~\ref{lemma_root_estimate} implies
\begin{corollary}\label{cor_norm_cons}
Let $M$ be a matrix with $\|MM^*-M^*M\|_{tr}\leq \e$ and $\|M\|_{op}\leq 1$.
Then there exists a normal matrix $N$ such that $\|M-N\|_{tr}\leq 36\e^{1/18}$ and
$\|N\|_{op}\leq 1$.
\end{corollary}
\begin{proof}
Let $M=UH$ with unitary $U$ and positive $H$. We have
$\|UH^2-H^2U\|_{tr}\leq\e$. So, by Theorem~\ref{th_2_self-adjoint_and_unitary}
we can find positive $A$ and unitary $V$ such that $\|H^2-A\|_{tr}\leq 30\e^{1/9}$,
$\|U-V\|_{tr}\leq 30\e^{1/9}$ and
$[H^2,A]=[V,A]=0$. By  Lemma~\ref{lemma_root_estimate} we have
$\|H-A^{1/2}\|\leq 6\e^{1/18}$
and $N=VA^{1/2}$ satisfies the Corollary.
\end{proof}

\section{Concluding remarks}
We see that the normalized Hilbert-Schmidt norm is more friendly for
almost-near questions for the commutator.  We think that it is
interesting to consider other relations. For example, if almost
solutions of
$$
U^k=V^{-1}UV
$$
are near solutions?

{\bf Acknowledgement} Some questions answered in the paper arise
during my talks at ``Nonstandard Analysis'' seminar at
Urbana-Champaign. I am thankful to E. Gordon, P. Leob and W. Henson
who was listening my messy talks and gave useful suggestions. I
think that I have had much more benefits from the talks then they
had.

The ``design'' of the
introduction is almost copied from \cite{Choi1}.

\end{document}